\documentclass[11pt,leqno]{article}
\usepackage[margin=1in]{geometry} 
\usepackage{amssymb,amsfonts,amsmath,bbm,mathrsfs,stmaryrd,mathtools}
\usepackage{xcolor}
\usepackage{url}

\usepackage{graphicx}

\usepackage{accents}

\usepackage{extarrows}

\usepackage[shortlabels]{enumitem}
\usepackage{tensor}

\usepackage{xr}
\usepackage[T1]{fontenc}
\usepackage[utf8]{inputenc}

\usepackage[colorlinks,
linkcolor=black!75!red,
citecolor=blue,
pdftitle={},
pdfproducer={pdfLaTeX},
pdfpagemode=None,
bookmarksopen=true,
bookmarksnumbered=true,
backref=page]{hyperref}

\usepackage{tikz}
\usetikzlibrary{arrows,calc,decorations.pathreplacing,decorations.markings,decorations.shapes,intersections,shapes.geometric,through,fit,shapes.symbols,positioning,decorations.pathmorphing}

\makeatletter
\newlength\zig@L
\newlength\zig@La
\newlength\zig@Lb

\newcommand{\xzigrightarrow}[2][]{%
  \mathrel{%
    \settowidth{\zig@La}{$\scriptstyle #2$}%
    \settowidth{\zig@Lb}{$\scriptstyle #1$}%
    \zig@L=\zig@La\relax
    \ifdim\zig@Lb>\zig@L \zig@L=\zig@Lb\fi
    \advance\zig@L by 2.2em\relax
    \tikz[baseline=-0.65ex]{%
      \draw[->,
            line cap=round,
            decorate,
            decoration={zigzag,segment length=4pt,amplitude=1.1pt}]%
        (0,0) -- (\zig@L,0)
        node[midway,above=2pt] {$\scriptstyle #2$}%
        \if\relax\detokenize{#1}\relax\else
          node[midway,below=2pt] {$\scriptstyle #1$}%
        \fi
      ;
    }%
  }%
}
\makeatother

\makeatletter
\newcommand{\squigjoin}{1mu} 

\def\sqleft@{\sim}                    
\def\sqmid@{\sim\mkern-\squigjoin}    

\def\rightsquigarrowfill@{%
  \arrowfill@{\sqleft@}{\sqmid@}{\mkern-4mu\succ}%
}

\newcommand{\xrightsquigarrow}[2][]{%
  \ext@arrow 0359\rightsquigarrowfill@{#1}{#2}%
}
\makeatother

\makeatletter
\newcommand*\circled[1]{\tikz[baseline=(char.base)]{
    \node[shape=circle, draw, inner sep=0pt, 
    minimum height={\f@size},] (char) {\vphantom{WAH1g}#1};}}
\makeatother

\makeatletter
\DeclareRobustCommand\widecheck[1]{{\mathpalette\@widecheck{#1}}}
\def\@widecheck#1#2{%
    \setbox\z@\hbox{\m@th$#1#2$}%
    \setbox\tw@\hbox{\m@th$#1%
       \widehat{%
          \vrule\@width\z@\@height\ht\z@
          \vrule\@height\z@\@width\wd\z@}$}%
    \dp\tw@-\ht\z@
    \@tempdima\ht\z@ \advance\@tempdima2\ht\tw@ \divide\@tempdima\thr@@
    \setbox\tw@\hbox{%
       \raise\@tempdima\hbox{\scalebox{1}[-1]{\lower\@tempdima\box
\tw@}}}%
    {\ooalign{\box\tw@ \cr \box\z@}}}
\makeatother

\usepackage{braket}

\usepackage[amsmath,thmmarks,hyperref]{ntheorem}
\usepackage{cleveref}

\newcommand\nthalias[1]{\AddToHook{env/#1/begin}{\crefalias{lemma}{#1}}}

\nthalias{definition}
\nthalias{example}
\nthalias{examples}
\nthalias{remark}
\nthalias{remarks}
\nthalias{convention}
\nthalias{notation}
\nthalias{construction}
\nthalias{sketch}
\nthalias{theoremN}
\nthalias{propositionN}
\nthalias{corollaryN}
\nthalias{lemma}
\nthalias{proposition}
\nthalias{corollary}
\nthalias{theorem}
\nthalias{conjecture}
\nthalias{question}
\nthalias{assumption}

\creflabelformat{enumi}{#2#1#3}

\crefname{section}{Section}{Sections}
\crefformat{section}{#2Section~#1#3} 
\Crefformat{section}{#2Section~#1#3} 

\crefname{subsection}{\S}{\S\S}
\AtBeginDocument{%
  \crefformat{subsection}{#2\S#1#3}%
  \Crefformat{subsection}{#2\S#1#3}%
}

\crefname{subsubsection}{\S}{\S\S}
\AtBeginDocument{%
  \crefformat{subsubsection}{#2\S#1#3}%
  \Crefformat{subsubsection}{#2\S#1#3}%
}

\theoremstyle{plain}

\newtheorem{lemma}{Lemma}[section]
\newtheorem{proposition}[lemma]{Proposition}
\newtheorem{corollary}[lemma]{Corollary}
\newtheorem{theorem}[lemma]{Theorem}

\theoremstyle{plain}
\theoremnumbering{Alph}

\theoremstyle{plain}
\theorembodyfont{\upshape}
\theoremsymbol{\ensuremath{\blacklozenge}}

\newtheorem{remark}[lemma]{Remark}
\newtheorem{remarks}[lemma]{Remarks}

\crefname{definition}{definition}{definitions}
\crefformat{definition}{#2definition~#1#3} 
\Crefformat{definition}{#2Definition~#1#3} 

\crefname{ex}{example}{examples}
\crefformat{example}{#2example~#1#3} 
\Crefformat{example}{#2Example~#1#3} 

\crefname{exs}{example}{examples}
\crefformat{examples}{#2example~#1#3} 
\Crefformat{examples}{#2Example~#1#3} 

\crefname{remark}{remark}{remarks}
\crefformat{remark}{#2remark~#1#3} 
\Crefformat{remark}{#2Remark~#1#3} 

\crefname{remarks}{remark}{remarks}
\crefformat{remarks}{#2remark~#1#3} 
\Crefformat{remarks}{#2Remark~#1#3} 

\crefname{convention}{convention}{conventions}
\crefformat{convention}{#2convention~#1#3} 
\Crefformat{convention}{#2Convention~#1#3} 

\crefname{notation}{notation}{notations}
\crefformat{notation}{#2notation~#1#3} 
\Crefformat{notation}{#2Notation~#1#3} 

\crefname{table}{table}{tables}
\crefformat{table}{#2table~#1#3} 
\Crefformat{table}{#2Table~#1#3}

\crefname{lemma}{lemma}{lemmas}
\crefformat{lemma}{#2lemma~#1#3} 
\Crefformat{lemma}{#2Lemma~#1#3} 

\crefname{proposition}{proposition}{propositions}
\crefformat{proposition}{#2proposition~#1#3} 
\Crefformat{proposition}{#2Proposition~#1#3} 

\crefname{propositionN}{proposition}{propositions}
\crefformat{propositionN}{#2proposition~#1#3} 
\Crefformat{propositionN}{#2Proposition~#1#3} 

\crefname{corollary}{corollary}{corollaries}
\crefformat{corollary}{#2corollary~#1#3} 
\Crefformat{corollary}{#2Corollary~#1#3} 

\crefname{corollaryN}{corollary}{corollaries}
\crefformat{corollaryN}{#2corollary~#1#3} 
\Crefformat{corollaryN}{#2Corollary~#1#3} 

\crefname{theorem}{theorem}{theorems}
\crefformat{theorem}{#2theorem~#1#3} 
\Crefformat{theorem}{#2Theorem~#1#3} 

\crefname{theoremN}{theorem}{theorems}
\crefformat{theoremN}{#2theorem~#1#3} 
\Crefformat{theoremN}{#2Theorem~#1#3} 

\crefname{enumi}{}{}
\crefformat{enumi}{#2#1#3}
\Crefformat{enumi}{#2#1#3}

\crefname{assumption}{assumption}{Assumptions}
\crefformat{assumption}{#2assumption~#1#3} 
\Crefformat{assumption}{#2Assumption~#1#3} 

\crefname{construction}{construction}{Constructions}
\crefformat{construction}{#2construction~#1#3} 
\Crefformat{construction}{#2Construction~#1#3} 

\crefname{sketch}{sketch}{Sketches}
\crefformat{sketch}{#2sketch~#1#3} 
\Crefformat{sketch}{#2Sketch~#1#3} 

\crefname{question}{question}{Questions}
\crefformat{question}{#2question~#1#3} 
\Crefformat{question}{#2Question~#1#3} 

\crefname{equation}{}{}
\crefformat{equation}{(#2#1#3)} 
\Crefformat{equation}{(#2#1#3)}

\numberwithin{equation}{section}

\theoremstyle{nonumberplain}
\theoremsymbol{\ensuremath{\blacksquare}}

\newtheorem{proof}{Proof}
\newcommand\pf[1]{\newtheorem{#1}{Proof of \Cref{#1}}}

\newcommand\bG{{\mathbb G}}
\newcommand\bH{{\mathbb H}}

\newcommand\bM{{\mathbb M}}

\newcommand\bR{{\mathbb R}}
\newcommand\bS{{\mathbb S}}

\newcommand\cC{{\mathcal C}}

\newcommand\cH{{\mathcal H}}

\newcommand\cM{{\mathcal M}}

\newcommand\cP{{\mathcal P}}

\newcommand\fg{{\mathfrak g}}

\DeclareMathOperator{\Ad}{Ad}
\DeclareMathOperator{\ad}{ad}

\DeclareMathOperator{\End}{\mathrm{End}}
\DeclareMathOperator{\Hom}{\mathrm{Hom}}

\DeclareMathOperator{\Aut}{\mathrm{Aut}}

\DeclareMathOperator{\im}{\mathrm{im}}

\newcommand{\qedhere}{\mbox{}\hfill\ensuremath{\blacksquare}}

\newcommand{\comment}[1]{}

\newcommand{\xrightarrowdbl}[2][]{%
  \xrightarrow[#1]{#2}\mathrel{\mkern-14mu}\rightarrow
}

\title{Controlled Riemannian rigidity}
\author{Alexandru Chirvasitu}

\begin{document}

\date{}

\newcommand{\Addresses}{{
  \bigskip
  \footnotesize

  \textsc{Department of Mathematics, University at Buffalo}
  \par\nopagebreak
  \textsc{Buffalo, NY 14260-2900, USA}  
  \par\nopagebreak
  \textit{E-mail address}: \texttt{achirvas@buffalo.edu}

}}

\maketitle

\begin{abstract}
  Given a compact Lie group acting smoothly and freely on a compact connected manifold, we prove that the space of invariant Riemannian metrics whose symmetry group is precisely the group in question is open dense in that of all invariant metrics provided the orbit space is at least 2-dimensional. This strengthens the author's prior results to their optimal dimension bound and recovers Ebin's non-equivariant analogue to the effect that trivial-symmetry Riemannian metrics constitute an open dense set. In the process, the discussion also addresses a number of apparent gaps in said sources, pertaining to the rigidity results in question. 
\end{abstract}

\noindent \emph{Key words:
  Riemannian manifold;
  compact Lie group;
  geodesic;
  isometry group;
  isotropy group;
  orbit;
  polar coordinates;
  principal action
}

\vspace{.5cm}

\noindent{MSC 2020: 58D17; 57S15; 53B20; 53C17; 53C07; 58A30; 53C22; 22D12
  
}


\section*{Introduction}

The running theme throughout the sequel is that of \emph{generic rigidity}, in the context of compact Lie groups acting isometrically on compact Riemannian manifolds. Specifically, having fixed a smooth compact-Lie action $M\circlearrowleft \bG$, the expectation might be that within the Fr\'echet manifold $\cM_{\bG}(M)$ of $\bG$-invariant Riemannian metrics, those whose isometries consist \emph{precisely} of $\bG$ constitute a ``large'' set (\emph{residual} \cite[\S 9]{oxt_meas-cat_2e_1980} or, better perhaps, open dense).

The pattern is a familiar one in the literature on realizing prescribed symmetry groups: it is effectively the strategy adopted in \cite[Theorem]{sab-inf}, say, of constructing graphs with a given automorphism group by breaking an initially larger symmetry group, the same goes for specifying polytope symmetries in \cite[Theorem 1]{MR4239243}, and there is no shortage of other examples. One of the note's main results, very much in this same spirit and direction, is as follows. 

\begin{theorem}\label{th:rig.lg.dim}
  Let $M\circlearrowleft \bG$ be a free $C^{\infty}$ action by a compact Lie group on a compact connected $C^{\infty}$ manifold. If $e:=\dim M/\bG\ge 2$ then the subset
  \begin{equation*}
    \left\{
      g\in \cM_{\bG}
      \ :\
      \Aut(g)=\bG
    \right\}
    \subseteq
    \cM_{\bG}
  \end{equation*}
  is open dense. 
\end{theorem}

This improves on \cite[Theorem 4.3]{zbMATH07342206} (which requires $e\ge \dim \bG+1$), in addition to repairing both a gap in that earlier result's proof and omissions in its statement, as elaborated in \hyperref[res:how.brkn]{Remarks~\ref*{res:how.brkn}}. The result hinges on first proving the genericity (in the sense of Fr\'echet dense openness again) for \emph{vertical} metrics: those whose symmetry group leaves all orbits invariant. This is achieved in the broader context of \emph{principal} actions (i.e. with only principal orbits \cite[\S IV.3]{bred_cpct-transf}). On the one hand, openness is automatic via

\begin{theorem}\label{th:vert.open}
  Given a principal action $M\circlearrowleft \bG$ of a compact Lie group on a compact manifold, the condition that a $\bG$-invariant Riemannian metric on $M$ be vertical is open in $\cM_{\bG}(M)$. 
\end{theorem}

On the other, density obtains for optimal dimension bounds:

\begin{theorem}\label{th:eqvr.vert.sharp.dim.2}
  Let $M\circlearrowleft \bG$ be a principal action of a compact Lie group on a compact smooth manifold. If 
  \begin{equation*}
    \forall\left(\text{component }M_i\le M\right)
    \forall\left(p\in M_i\right)
    \left(
      \dim M_i\ge 2+\dim p\bG
    \right)
  \end{equation*}
  then the space $\cM^{\uparrow}_{\bG}(M)\subseteq \cM_{\bG}(M)$ of invariant vertical metrics is open dense. 
\end{theorem}

Leveraging the verticality package provided by \Cref{th:vert.open,th:eqvr.vert.sharp.dim.2} into \Cref{th:rig.lg.dim} is effected by means of results (\Cref{th:gen.onto.curv,th:gen.t.tens,th:r.and.t}) concerning the genericity of the space of maximal-rank Ehresmann connections and the freedom to prescribe O'Neill tensors associated to $\bG$-invariant metrics, possibly of some independent interest.

\subsection*{Acknowledgments}

The original motivation for the investigation was an audit run on \cite{zbMATH07342206} by OpenAI's GPT-6, uncovering the gap marring \cite[Claim 4.7]{zbMATH07342206}. The present paper's text, proof strategies and arguments are the author's, with technical audits for correctness performed by the same LLM. The author assumes full responsibility for the resulting work.


\section{Generic verticality and rigidity for compact-Lie-invariant Riemannian metrics}\label{se:gen.vert.rig}

We write $\cM(M)$ for the space of (smooth, i.e. $C^{\infty}$) Riemannian metrics on a $C^{\infty}$ manifold $M$, and $\cM_{\bG}(M)$ for that of $\bG$-invariant metrics for a fixed smooth $\bG$-action $M\circlearrowleft \bG$ by a (typically compact) Lie group $\bG$. 

Recall (e.g. \cite[post Proposition I.2.5]{audin}) that a compact-Lie-group action $M\circlearrowleft \bG$ on a manifold with connected $M/\bG$ is \emph{principal} if its orbits all are, i.e. are among the orbits whose \emph{type} (i.e. isotropy-group conjugacy class \cite[post Corollary I.4.4]{bred_cpct-transf}) is the single one encountered open densely.

\pf{th:vert.open}
\begin{th:vert.open}
  Consider a vertical $g\in \cM_{\bG}(M)$; coupled with $\bG$-invariance, verticality means that the orbits of $\Aut(g)$ are precisely the $\bG$-orbits. For $\cM_G(M)\ni g'\sim g$ (sufficiently $C^{\infty}$-close), the corresponding $\Aut(g')$-orbits
  \begin{itemize}[wide]
  \item contain $\bG$-orbits (again by the assumed $\bG$-invariance);
    
  \item and are diffeomorphic to compact smooth submanifolds of the $\Aut(g)$-orbits by Ebin's upper semicontinuity theorem \cite[Theorem 8.1]{zbMATH03327104} (also \cite[Proposition 4.1]{MR4344307}):
    \begin{equation*}
      \Ad_{\sigma} \Aut(g')
      :=
      \sigma \Aut(g')\sigma^{-1}
      \le
      \Aut(g)
      ,\quad
      1\sim \sigma\in \Aut_{\cC^{\infty}}(M).
    \end{equation*}
  \end{itemize}
  An arbitrary $\Aut(g')$-orbit $O$ will thus be a disjoint union of $k'$ closed (compact, connected, boundary-less) $C^{\infty}$ $q'$-manifolds for $q'\ge q=\dim p\bG$ (arbitrary $p$). An embedding $O\subseteq p\bG$ must identify each of those components with one of the $k$ components of $p\bG$ by compactness and \emph{domain invariance} \cite[Corollary 2B.4]{hatch_at}, so that in particular $q=q'$ and $k=k'$ and $\Aut(g')$-orbits must be $\bG$-orbits.
\end{th:vert.open}

Coupled with \cite[Theorem 3.2]{zbMATH07342206}, \Cref{th:vert.open} proves

\begin{corollary}\label{cor:vert.op.dns}
  If $M\circlearrowleft \bG$ is a principal action of a compact Lie group on a compact smooth manifold whose components all have dimensions $\ge 3+\dim \bG$ then the space
  \begin{equation*}
    \cM_{\bG}^{\uparrow}(M)
    :=
    \left\{
      g\in \cM_{\bG}(M)
      \ :\
      \text{$\Aut(g)$-orbits are (contained in) $\bG$-orbits}
    \right\}
  \end{equation*}
  of \emph{vertical} $\bG$-invariant Riemannian metrics is open dense in $\cM_{\bG}(M)$.  \qedhere
\end{corollary}

\begin{remark}\label{re:bd23}
  The sharp dimension bound for $\bG=\{1\}$ of \Cref{cor:vert.op.dns} is in fact $\dim M=\dim M-\dim \bG=2$ (rather than the statement's stipulated $3$), by \cite[Proposition 8.3]{zbMATH03327104}. \cite[Remark 3.3]{zbMATH07342206} further notes that a gap of 1 is certainly insufficient, leaving open the possibility that the optimal value might be 2. \Cref{th:eqvr.vert.sharp.dim.2} confirms this, improving on \cite[Theorem 3.2]{zbMATH07342206} in that respect. 
\end{remark}

\pf{th:eqvr.vert.sharp.dim.2}
\begin{th:eqvr.vert.sharp.dim.2}
  There is no loss in assuming the action faithful if so desired (for one can always substitute $\bG/\left(\ker~\text{action}\right)$ for $\bG$). Moreover, $M/\bG$ may as well be connected: small perturbations will equip the quotient with a metric giving its components distinct geodesic diameters.

  Openness is already settled by \Cref{th:vert.open}, so the core issue is density. We will adapt, to that purpose, the argument employed in proving \cite[Proposition 8.3]{zbMATH03327104}; the proof distinguishes a number of stages.
  
  \begin{enumerate}[(I),wide]
  \item\label{item:th:eqvr.vert.sharp.dim.2:pf.glob.inv.orb}\textbf{: a global invariant orbit $p\bG$.} One can be assumed via the gadget operative in \cite[Lemma 3.8]{zbMATH07342206}: any given $\bG$-invariant metric can be perturbed (equivariantly) as slightly as one desires to another, say $g$, so as to ensure maximization of the scalar-curvature absolute value along a single orbit $p\bG$ for some $p$. That orbit, then, will be invariant under the symmetry group $J_g:=\Aut(g)$.

    Note that the cited \cite[Lemma 3.8]{zbMATH07342206} makes the stronger assumption that $M$ have dimension dominating $2+\dim \bG$, but the argument in fact only requires the bound
    \begin{equation*}
      e:=
      \inf_{\substack{\text{component $M_i$}\\p\in M_i}}
      \left(\dim M_i-\dim p\bG\right)
      \ge 2
    \end{equation*}
    (i.e. the present hypothesis).
    
  \item\label{item:th:eqvr.vert.sharp.dim.2:pf.ind}\textbf{: induction on dimension.} Assume $e\ge 3$, with the statement settled for all smaller values $\ge 2$ thereof. Having first perturbed the original metric so as to ensure the existence of a unique point (there denoted by $q$) where the maximal eigenvalue of the Ricci tensor is achieved, the proof of \cite[Proposition 8.3]{zbMATH03327104} then proceeds to further perturb the metric ``around'' a small-radius geodesic sphere centered at $q$. \cite[\S 3.1]{zbMATH07342206} flags a potential gap in that proof (or at any rate, at least in its presentation): there is no indication as to how that perturbation is to be performed, and apparently no guarantee that the procedure will preserve geodesic spheres (so that the isometry group of the new metric preserves them). To mitigate that issue, proceed as follows:
    \begin{itemize}[wide]
    \item perturb the current metric
      \begin{equation*}
        g=dr^2+g_r
        ,\quad
        r:=\text{radial coordinate away from $p\in p\bG$},
      \end{equation*}
      assumed to leave invariant the orbit $p\bG$ produced by \Cref{item:th:eqvr.vert.sharp.dim.2:pf.glob.inv.orb} above, only along the horizontal manifold $H_{g,p}$ consisting of geodesics originating at $p$;
      
    \item with the perturbation only along the angular coordinates in the \emph{polar Fermi coordinate system} \cite[Example 6.46]{lee_riem_2e_2018} around the orbit $p\bG$ (i.e. $g_r$ is perturbed):
      \begin{equation*}
        g'_r-g_r=\varphi(r)\varepsilon
        ,\quad
        \bR_{\ge 0}
        \xrightarrow[\quad\text{smooth}\quad]{\quad\varphi\quad}
        [0,1]
      \end{equation*}
      for a bump function $\varphi$ supported on a small interval around some fixed small $r_0$. 
    \end{itemize}
    With these caveats, the geodesic spheres
    \begin{equation*}
      H_{g,p,r}:=
      \left\{x\in H_{g,p}\ :\ d_g(p,x)=r\quad\left(\text{geodesic distance}\right)\right\}
    \end{equation*}
    remain so. The perturbation can be made so as to produce a rigid metric along the (at least 2-dimensional) sphere $H_{g,p,r_0}$. This then ensures that the isotropy group $J_{g,p}$ (of $p$, with respect to the newly obtained $g$) fixes $\left(T_p p\bG\right)^{\perp}$ pointwise, which suffices to complete the target inductive step:
    \begin{itemize}[wide]
    \item a fixed global orbit $p\bG$;
    \item coupled with the trivial action of $J_{g,p}$ on the horizontal tangent space $\left(T_p p\bG\right)^{\perp}$;
    \item and with the assumed connectedness of $M/\bG$
    \end{itemize}
    translates to verticality. 
    
  \item\label{item:th:eqvr.vert.sharp.dim.2:pf.base}\textbf{: the base case $e=2$.} This will again be an angular perturbation, this time around the \emph{circle} $H_{g,p,r_0}$. The orthogonal group $O\left(\left(T_p p\bG\right)^{\perp}\right)$ will operate on that circle for it does on the entire disk $H_{g,p,\le r_0}$, and it will suffice to perturb the function
    \begin{equation*}
      H_{g,p,r_0}
      \ni
      q=q(\theta,r_0)
      \xmapsto{\quad}
      a(\theta)
      ,\quad
      g_{r_0}=a(\theta)d\theta^2
      \quad\left(\text{polar coordinates}\right)
    \end{equation*}
    arbitrarily slightly to a smooth function $a(1+\varepsilon)$ with trivial isotropy in $O\left(\left(T_p p\bG\right)^{\perp}\right)$. Said function will be $q(\theta,r_0)\mapsto a(\theta)g'_{r_0}(X_q,X_q)$ with $g$-unit $X_q\in T_q H_{g,p,r_0}$ for the angular component $g'_{r_0}$ of the perturbed metric. The desired trivial-isotropy deformation is indeed possible by \Cref{pr:gen.triv.iotrop} below, hence the conclusion.  \qedhere
  \end{enumerate}
\end{th:eqvr.vert.sharp.dim.2}

\begin{remarks}\label{res:rec.rig.triv.g}
  \begin{enumerate}[(1),wide]
  \item Note incidentally that \Cref{th:eqvr.vert.sharp.dim.2} recovers \cite[Proposition 8.3]{zbMATH03327104} for trivial $\bG$ (so constitutes an equivariant generalization thereof): the proof implements a strategy very close in spirit to what the earlier result presumably intended, but does not rely on that result itself.

  \item \cite[Lemma 3.8]{zbMATH07342206}, on which the proof of \Cref{th:eqvr.vert.sharp.dim.2} is reliant, gives a somewhat misleading account of one of its claims (though in a manner apparently not directly relevant to the uses the lemma is put to above). The distinction between the earlier result and \Cref{le:le38.uncond} is that the former makes the $U$-localization aspect of the deformation conditional on an unnecessary constraint, and the proposed proof provides insufficient justification under that constraint. 
  \end{enumerate}
\end{remarks}

\begin{lemma}\label{le:le38.uncond}
  Let $M\circlearrowleft \bG$ be a smooth action of a compact Lie group on a compact $C^{\infty}$ manifold with connected components of dimension $\ge 2+\dim \bG$, and $g\in \cM_{\bG}(M)$. There is a $p\in M$ such that
  \begin{itemize}[wide]
  \item arbitrarily small perturbations $g'\in \cM_{\bG}(M)$ of $g$ exist, coincident with it outside arbitrarily small neighborhoods $U=U\bG\supset p\bG$;
  \item maximizing scalar-curvature absolute values on a unique orbit in $U$;
  \item and consequently with $\Aut(g')$ preserving said orbit. 
  \end{itemize}
\end{lemma}
\begin{proof}
  Select $p\in M$ maximizing the absolute scalar curvature $\left|R_{g}\right|$ and proceed as in the original proof's iterative procedure, with the added constraint of working only within $U$ (leaving $g$ unaltered outside that neighborhood): perturb $g$ conformally and iteratively through $g_{n}\mapsto \varphi_{n}^{-2}g_n$ for positive $\bG$-invariant smooth functions $\varphi_n$ on $M$, chosen judiciously so as to ensure a limit exhibiting the requisite curvature uniqueness in $U$. That the sought-after perturbations are achievable follows from the scalar-curvature variation formula \cite[Theorem 1.159(f)]{bss_einst_2008} under conformal alteration:
  \begin{equation*}
    R_{e^{-2f}g}
    =
    e^{2f}
    \left(R - 2\left(\dim M-1\right)\Delta f-(\dim M-1)(\dim M-2)\left|\mathrm{d} f\right|^2\right)
  \end{equation*}
  also invoked in the course  of the proof of \cite[Proposition 3.6(b)]{zbMATH07342206}.
\end{proof}

\begin{proposition}\label{pr:gen.triv.iotrop}
  Given a smooth action $M\circlearrowleft \bG:\alpha$ by a compact Lie group on a smooth metrizable manifold, the space of smooth functions on $M$ whose isotropy group is precisely $\ker\alpha$ is $C^{\infty}$-open dense. 
\end{proposition}
\begin{proof}
  It will be harmless to assume $\alpha$ faithful, the target being \emph{trivial} isotropy. Openness obtains quite broadly:
  \begin{itemize}[wide]
  \item on the one hand, for any continuous compact-group action $X\times \bG\to X$ the isotropy-group map
    \begin{equation*}
      X\ni x
      \xmapsto{\quad}
      \bG_x:=\left\{s\in \bG\ :\ xs=x\right\}
      \in
      \left\{
        \text{closed subgroups }\le \bG
      \right\}
    \end{equation*}
    is easily proven \emph{upper semicontinuous} \cite[Definition 7.1.1]{kt_corresp}. This suffices to conclude that the trivial-isotropy $x\in X$ constitute an open set when $\bG$ is Lie, having an identity neighborhood which contains no non-trivial subgroups \cite[Proposition 2.17]{hm_pro-lie-bk}.

  \item On the other hand, our continuous $M\circlearrowleft \bG$ certainly induces a \emph{separately} continuous action $C^{\infty}(M)\circlearrowleft \bG$, hence \cite[\S VIII.2.2, Proposition 1]{bourb_int_en_7-9} also jointly continuous on the \emph{barreled} space $C^{\infty}(M)$. The latter condition obtains under the metrizability hypothesis:
    \begin{itemize}[wide]
    \item when $M$ is connected (hence $\sigma$-compact \cite[Theorem 2.1]{gauld_non-met-mfld}) by \cite[Corollary 6.9]{MR4398224} because $C^{\infty}$ is then even \emph{Fr\'echet} \cite[Corollary 1 to Proposition 33.2]{trev_tvs};

    \item so for arbitrary metrizable $M$ by barreled-space permanence under products \cite[\S 27.1(5)]{k_tvs-1}, given the decomposition
      \begin{equation*}
        C^{\infty}(M)
        \cong
        \prod_{\text{components $M_i$}}C^{\infty}(M_i).
      \end{equation*}      
    \end{itemize}
  \end{itemize}
  It thus remains to settle the density claim. It, moreover, reduces to homogeneous manifolds $\bH\backslash \bG \circlearrowleft \bG$: assuming the latter variant settled, one can perturb any given $f\in C^{\infty}(M)$ slightly along each of finitely many orbits
  \begin{equation*}
    p_i \bG\cong \bH_i\backslash\bG
    ,\quad
    1\le i\le k
    ,\quad
    \bigcap_{i=1}^k \bigcap_{s\in \bG} \Ad_s\bH_i=\ker\alpha=\{1\},
  \end{equation*}
  the latter trivial intersection being achievable because descending chains of closed subgroups in compact Lie groups eventually stabilize. The perturbation can then be effected globally by Whitney extension \cite[Theorem I.4.1]{malgr_ideals}. For further convenience:
  \begin{itemize}[wide]
  \item $M=\bH\backslash \bG$ may as well be connected, for a function can always be slightly shifted along individual components so as to ensure it achieves distinct maxima thereon;
    
  \item $\dim \bH\backslash \bG>0$, for there is nothing to prove for (now connected) singletons $\bH\backslash \bG$;

  \item and the action can, recall, be assumed faithful: substitute $\bG/\left(\text{kernel}\right)$ for $\bG$.     
  \end{itemize}
  The remainder of the proof will consist of a recursive process very much reminiscent of \Cref{th:eqvr.vert.sharp.dim.2}'s induction (and \cite[Proposition 8.3]{zbMATH03327104}'s proof). Specifically, having selected a target $f\in C^{\infty}$ to deform arbitrarily slightly into a trivial-isotropy function, proceed as follows.

  \begin{itemize}[wide]
  \item Fix a $\bG$-invariant Riemannian metric $g$ on $M$, assumed available for the duration of the proof.
    
  \item On any smooth compact boundary-less manifold, the smooth functions achieving unique \emph{non-degenerate} maxima constitute an open dense subset of $C^{\infty}(M)$ \cite[Corollary 3.6]{morgan_morse_2020}. $f$ can thus be perturbed so as to ensure such a unique maximum at $x_0$.

  \item Further deform $f$ in a sufficiently small neighborhood of the radius-$r_0$ $x_0$-centered geodesic sphere $S_0:=M_{x_0,r_0}$ contained in a normal coordinate neighborhood around $x_0$ so as to ensure a unique global maximum at $x_1\in S_0$.

  \item We henceforth work entirely on $S_0$. For a small angle $0<\theta_1<\frac{\pi}{2}$ consider the smooth $(n-2)$-sphere $S_{1}\subseteq S_0$ cut out by the $x_0$-based geodesic rays making an angle of $\theta_1$ with the ray $x_0 x_1$ (indeed a sphere $C^{\infty}$-wise, because a level set in a polar-coordinate system around $x_0$); perturb $f$ so as to have it achieve a unique global maximum at $x_2\in S_1$. 
    
  \item Continue the procedure, up to and including $x_{n:=\dim M}$.

    The step producing $x_{k+1}$ will have available $x_{1\le \bullet\le k}$ fixed by the isotropy group $\bG_f$. The locus in $S_0$ cut out by the $x_0$-based geodesic rays making respective angles of $0<\theta_j<\frac {\pi}2$, $1\le j\le k$ (which we have the freedom to choose) with
    \begin{equation*}
      \mathrm{span}\left\{
        \text{$x_0$-velocity along }x_0 x_{\ell}
        \ :\
        1\le \ell\le j
      \right\}
      \le
      T_{x_0}M
    \end{equation*}
    is a disjoint union of $\bS^{n-k-1}$ if each angle $\theta_{\bullet}$ is chosen to be sufficiently small, having fixed the preceding ones. All such are individually invariant under $\bG_f$, for the locus cut out by fixing $\theta_{j}$ alone is a $\bS^{j-1}\times \bS^{n-j-1}$ with the point-wise fixed span of $x_0 x_{1\le \bullet\le j}$ as one of the $\bS^{j-1}$ slices. One can now select one of the $\bG_f$-invariant $S_k\cong \bS^{n-k-1}$ and maximize $f$ uniquely along it at $x_{k+1}$ by a small perturbation.
    
    \item Subsequent perturbations can always be chosen so small as to respect maxima fixed earlier and the angles $\theta_k$ can always be chosen sufficiently small to ensure the non-emptiness of the $S_k$. 
  \end{itemize}
  
  The isotropy group of the resulting $f$ in $\bG\subseteq \Aut(g)$ will fix all $x_i$, so must be trivial by the assumed connectedness of $M$: it fixes $x_0$ and operates trivially on $T_{x_0}M$, the $x_0$-velocities of $x_0 x_k$, $1\le k\le n$ being by construction linearly independent. 
\end{proof}

\begin{remark}\label{re:glb.lin.indep}
  Note the family resemblance between the ray linear independence in the proof of \Cref{pr:gen.triv.iotrop} and the analogue in that of \cite[Proposition 8.3]{zbMATH03327104}, where the counterparts to the $S_k$ are intersections of geodesic spheres (and the requisite linear independence is left to the reader).
\end{remark}

\pf{th:rig.lg.dim}
\begin{th:rig.lg.dim}
  The argument is structured so as to emphasize that \Cref{th:gen.onto.curv} already provides most of what is needed independently of \Cref{th:gen.t.tens}, with only a slightly sub-optimal bound for $e$.

  \begin{enumerate}[label={},wide]
  \item\textbf{($e\ge \min\left(3,\max\left(2,1+\dim \bG\right)\right)$)} Openness is noted in \cite[Theorem 4.3]{zbMATH07342206} (and covered by that result's valid portion), so the substance is density. We may also assume throughout the proof, per \Cref{th:eqvr.vert.sharp.dim.2}, that we work with vertical metrics only (for those constitute an open set, so small perturbations will not affect verticality). 
    
    A $\bG$-invariant Riemannian structure consists of the following pieces of data (cf. \cite[\S 9.15]{bss_einst_2008}):
    \begin{itemize}[wide]
    \item a Riemannian metric on $M/\bG$, base space of the principal bundle \cite[Theorem II.5.8]{bred_cpct-transf} $M\xrightarrowdbl{\pi}M/\bG$ attached to the free action;

    \item a smooth family of right-$\bG$-invariant metrics on the fibers of that principal bundle (i.e. orbits);

    \item and the \emph{horizontal distribution} $\cH_g$ \cite[\S 9.7]{bss_einst_2008}
      \begin{equation*}
        \bigsqcup_{p\in M}\left(T_pp\bG\right)^{\perp}
        \xrightarrowdbl{\quad}
        M.
      \end{equation*}
    \end{itemize}
    The latter, in particular, is nothing but a \emph{connection} \cite[Definition 9.39]{bss_einst_2008} on the principal $\bG$-bundle $\pi$.

    Having fixed a $g\in \cM^{\uparrow}_{\bG}$ (to be perturbed slightly into one with automorphism group $\bG$), observe that
    \begin{enumerate}[(a),wide]
    \item\label{item:th:rig.lg.dim:pf.inv.fib} fibers are in any case invariant under $J_g:=\Aut(g)$ by the assumed vertical character of $g$;

    \item\label{item:th:rig.lg.dim:pf.conn.m} the connectedness of $M$ reduces the goal to showing that for any one $p\in M$ the isotropy subgroup $J_{g,p}\le J_g$ is trivial;

    \item\label{item:th:rig.lg.dim:pf.vert} and, given the triviality of the action of $J_{g,p}$ on the manifold $H_{g,p}$ spanned by $p$-based geodesic arcs sufficiently short to avoid singularities (once again a consequence of verticality), the desired conclusion in fact reduces to proving the action $J_{g,p}\circlearrowright T_p p\bG$ trivial.
    \end{enumerate}
    More can be said concerning \Cref{item:th:rig.lg.dim:pf.vert}, which can be iterated: $J_{g,p}$ operates trivially on the horizontal tangent space $\cH_{g,q}$ at any point $q$ accessible from $p$ through broken horizontal curves, hence also on iterated Lie brackets of local horizontal vector fields around such points.

    One consequence of the above discussion: per \emph{Ambrose-Singer} \cite[Theorem D.2.1]{mntgmr_subriem_2002}, $J_{g,p}$ fixes pointwise the image of the \emph{curvature} (\cite[Definition 4.1.1]{mntgmr_subriem_2002}, \cite[\S 3]{MR702806})
    \begin{equation*}
      R^{\cH}\in \Gamma \left(\Lambda^2 T^* \left(M/\bG\right)\otimes \ad (M)\right)
      ,\quad
      \begin{aligned}
        \ad(M)
        &:=M\times_{\bG}\left(\fg:=Lie(\bG)\right)\\
        &:=\left(M\times \fg\right)/\bG
          \quad\left(\text{adjoint action $\bG\circlearrowright \fg$}\right)
      \end{aligned}
    \end{equation*}
    attached to the connection (horizontal distribution) $\cH_g$. To conclude:
    \begin{enumerate}[(I),wide]
    \item On the one hand, for $\dim \bG\le 1$ the assumption that $e\ge 2$ ensures via \Cref{th:gen.onto.curv} that arbitrarily slight $C^{\infty}$ perturbations of $g$ can render its curvature surjective at any one $p\bG\in M/\bG$.

    \item For $\dim \bG=2$ we have $e\ge 3$ and the curvature can again be chosen onto under small perturbations.
      
    \item As to $\dim\bG\ge 3$, with $e\ge 3$ again, a slight elaboration on the argument just given will do. The same \Cref{th:gen.onto.curv} allows for small deformations with $(\ge 3)$-dimensional $\im R_{p\bG}^{\cH}\le \fg$ arbitrary in a sufficiently small Grassmannian neighborhood of some $Z\le \fg$ containing the curvature image being deformed. The conclusion follows from the observation that generic $(\ge 3)$-dimensional subspaces of compact Lie algebras generate the latter under iterated Lie brackets and passing to Lie algebras of integral-group closures:
      \begin{itemize}[wide]
      \item a generic $v\in \fg$ generates a maximal torus topologically;

      \item two generic vectors in a semisimple Lie algebra generate the latter \cite[Theorem 6 and its proof]{zbMATH03070481};

      \item and finally, arbitrary compact Lie algebras decompose \cite[Theorems 6.15 and 6.18]{hm5} as sums of the latter two types. 
      \end{itemize}
    \end{enumerate}
    
  \item\textbf{($e=2$, $\dim \bG\ge 2$)} \Cref{th:r.and.t} supplies a $J_{g,p}$-invariant $0\ne v_0\in T_p p\bG$ (any $v_0$ in the image of the curvature will do) with 
    \begin{equation*}
      v_{j}
      ,\quad
      0\le j\le \dim\bG -1
      ,\quad
      v_{j+1}:= T(v_j,w)
      \quad\left(\text{some fixed $0\ne w\in\left(T_pp\bG\right)^{\perp}$ }\right)
    \end{equation*}
    linearly independent. 
  \end{enumerate}
  This completes the proof. 
\end{th:rig.lg.dim}

\begin{remarks}\label{res:how.brkn}
  \begin{enumerate}[(1),wide]
  \item One problem with \cite[Theorem 4.3]{zbMATH07342206} as appearing there (the less serious one) is its omission of any connectedness constraints. Assuming \Cref{th:rig.lg.dim}, the generic $\bG$-invariant metric on $M$ will have a power $\bG^{\left|\pi_0 M\right|}$ as its isotropy. That disconnected $M$ are easily accommodated for connected $\bG$ in this manner for the purpose of determining their generic automorphism group follows from the simple remark recorded in \Cref{le:disc.mfld}.
    
  \item More seriously, \cite[Claim 4.7]{zbMATH07342206}, to the effect that two $(\dim M/\bG)$-dimensional tangent spaces at some $q\in M$ can be chosen in general position, in fact proceeds to render the target spaces coincident.
  \end{enumerate}
\end{remarks}

\begin{lemma}\label{le:disc.mfld}
  Let $M\circlearrowleft \bG$ be a smooth action by a compact connected Lie group on a smooth compact manifold with $(>0)$-dimensional connected components $M_i$. If
  \begin{equation*}
    \forall i\quad:\quad
    \left\{g\in \cM_{\bG}(M_i)\ :\ \Aut(g)\cong \bG_i\right\}
    \quad
    \overset{\text{residual}}{\subseteq}
    \quad
    \cM_{\bG}(M_i)
  \end{equation*}
  for compact Lie groups $\bG_i$, then $\Aut(g)\cong \prod_i \bG_i$ over a residual set of $g\in \cM_{\bG}(M)$.
  
  The same goes, moreover, with ``open dense'' in place of ``residual'' throughout. 
\end{lemma}
\begin{proof}
  Indeed, slight metric rescaling will ensure distinct geodesic diameters for the components $M_i$, and hence their individual preservation by $\Aut(g)$ as an open dense condition. 
\end{proof}

\begin{theorem}\label{th:gen.onto.curv}
  Let $M\circlearrowleft \bG$ be a free $C^{\infty}$ action by a compact Lie group on a $C^{\infty}$ compact connected manifold, $p\in \bM$ and $e:=\dim M/\bG$.

  \begin{enumerate}[(1),wide]
  \item\label{item:th:gen.onto.curv:nbhd} For any
    \begin{equation*}
      \cH\in
      \mathrm{Conn}\left(M\circlearrowleft\bG\right)
      :=
      \left\{\text{connections on $M\circlearrowleft\bG$}\right\}
    \end{equation*}
    and any $C^{\infty}$ neighborhood $U\ni \cH$, the image $\im R^{\cH'}_{p\bG}\le \fg$, $\cH'\in U$ can be chosen arbitrarily within a sufficiently small Grassmannian neighborhood of any intermediate subspace
    \begin{equation*}
      \im R^{\cH}_{p\bG}
      \le
      Z
      \le
      \fg
      ,\quad
      \dim Z=
      \min\left(\dim\bG,\frac{e(e-1)}{2}\right).
    \end{equation*}
    
  \item\label{item:th:gen.onto.curv:many.conn} The set of connections
    \begin{equation*}
      \left\{
        \cH
        \ :\
        \dim\im R^{\cH}_{p\bG}
        =
        \min\left(\dim\bG,\frac{e(e-1)}{2}\right)
      \right\}
      \subseteq
      \mathrm{Conn}\left(M\circlearrowleft\bG\right)
    \end{equation*}
    is $C^{\infty}$-open dense. 
  \end{enumerate}
\end{theorem}
\begin{proof}
  That \Cref{item:th:gen.onto.curv:many.conn} follows from \Cref{item:th:gen.onto.curv:nbhd} is immediate, so the proof will dispatch the latter. Observe also that the target value is the theoretical maximum for the rank of a linear map $\Lambda^2 \bR^{e}\to \bR^{\dim \bG}$, so that openness is automatic in any case due to the continuity of
  \begin{equation*}
    \mathrm{Conn}\left(M\circlearrowleft\bG\right)
    \ni \cH
    \xmapsto{\quad}
    R^\cH_{p\bG}
    \in
    \Hom\left(\wedge^2 \bR^{e}, \bR^{\dim \bG}\right)
  \end{equation*}
  as a map of Fr\'echet manifolds and the lower semicontinuity of the rank function on the codomain.
  
  Having equipped $\mathrm{Conn}\left(M\circlearrowleft\bG\right)$ with its affine structure based on $\Gamma \left(T^* \left(M/\bG\right)\otimes \ad(M)\right)$ described in \cite[\S 3]{MR702806}, recall the resulting curvature variation formula for connection segments:
  \begin{equation*}
    \forall
    \left(
      \begin{aligned}
        \cH&\in \mathrm{Conn}\left(M\circlearrowleft\bG\right)\\
        \theta&\in \Gamma \left(T^* \left(M/\bG\right)\otimes \ad(M)\right)
      \end{aligned}
    \right)
    \left(
      R^{\cH+t\theta}
      =
      R^{\cH}
      +
      t d_{\cH} \theta
      +
      \frac{1}{2}
      t^2 [\theta,\theta]
    \right),
  \end{equation*}
  with
  \begin{equation*}
    \Gamma \left(\Lambda^{\bullet} T^* \left(M/\bG\right)\otimes \ad(M)\right)
    \xrightarrow{\quad d_{\cH}\quad}
    \Gamma \left(\Lambda^{\bullet+1} T^* \left(M/\bG\right)\otimes \ad(M)\right)
  \end{equation*}
  denoting the differential operator induced by covariant differentiation with respect to $\cH$ and $[\theta,\theta]$ the bracket on $\fg$-valued forms aggregating exterior products on the base and the Lie bracket on $\fg$. 

  \begin{enumerate}[(I),wide]
  \item\label{item:th:gen.onto.curv:pf.arb.diff}\textbf{: $\left(d_{\cH}\theta\right)_{p\bG}$ can be chosen arbitrarily.} Recall \cite[(4.2)]{MR702806} that
    \begin{equation*}
      d_{\cH}\theta(X,Y)
      =
      \nabla^{\cH}_X\left(\theta(Y)\right)
      -
      \nabla^{\cH}_Y\left(\theta(X)\right)
      -
      \theta([X,Y])
      ,\quad
      \begin{aligned}
        X,Y&\in \Gamma\left(T\left(M/\bG\right)\right)\\
        \nabla^{\cH}&:=\text{$\cH$-covariant derivative}.
      \end{aligned}
    \end{equation*}
    Allowing $X$ and $Y$ to range over $\partial_i:=\frac{\partial}{\partial x^i}$, $1\le i\le e$, the brackets $[X,Y]$ vanish. It remains to observe that the individual $\nabla^{\cH}_{\partial_i}\left(\theta(\partial_j)\right)$ can indeed be selected arbitrarily at the individual point $p\bG\in M/\bG$ as an application of, say, Whitney extension \cite[\S I.4, Theorem 4.1]{malgr_ideals}. 
    
  \item\label{item:th:gen.onto.curv:pf.end}\textbf{: Conclusion.} Item \Cref{item:th:gen.onto.curv:pf.arb.diff} above ensures that, for a fixed connection $\cH$, the smooth Fr\'echet-manifold map
    \begin{equation*}
      \Gamma\left(T^* \left(M/\bG\right)\otimes \ad(M)\right)
      \ni      
      \theta
      \xmapsto{\quad}
      R^{\cH+\theta}_{p\bG}
      \in
      \Hom\left(\Lambda^2 T_{p\bG}\left(M/\bG\right),\fg\right)
    \end{equation*}
    has surjective differential. The target being finite-dimensional, this suffices \cite[\S 1.56]{Schm23} to conclude that that map is a submersion locally around $\theta=0$ and its image $R^{\cH}_{p\bG}$. That map's openness, in particular, yields the desired conclusion. 
  \end{enumerate}
\end{proof}

For the purpose of lowering \Cref{th:rig.lg.dim}'s bound to the optimal $e\ge 2$, we need an analogue of \Cref{th:gen.onto.curv} above for the tensor $T$ attached in \cite[\S 9.C]{bss_einst_2008} (following \cite[\S 2]{MR200865}) to a $\bG$-invariant metric $g\in \cM_{\bG}(M)$: denoting by $\bullet^{\uparrow}$ and $\overrightarrow{\bullet}$ the operators respectively isolating the vertical/horizontal components of a vector field (perhaps with an additional ``$g$'' superscript for clarity),
\begin{equation*}
  T^{g}(X,Y)
  :=
  \left(
    \overrightarrow{T^g}(X,Y)
    :=
    \overrightarrow{\nabla_{X^{\uparrow}}Y^{\uparrow}}
  \right)
  +
  \left(
    T^{g\uparrow}(X,Y)    
    :=
    \left(\nabla_{X^{\uparrow}}\overrightarrow{Y}\right)^{\uparrow}
  \right)
  \in
  \cH\oplus \cH^{\perp}=TM. 
\end{equation*}
Dualizing and reversing arrow directions freely by means of the metric, the two components $\overrightarrow{T}$ and $T^{\uparrow}$ are, by \cite[(9.18c) and (9.18d)]{bss_einst_2008}, the images of a single section through the respective resulting bundle maps
\begin{equation*}
  S^2 T\ad(M)\otimes T^*\left(M/\bG\right)
  \ \xrightarrow{\quad}\ 
  T\left(M/\bG\right)\otimes T^{*\otimes 2}\ad(M)
  \quad\text{and}\quad
  \End\ad(M)\otimes T^*\left(M/\bG\right).
\end{equation*}

\begin{remark}\label{re:many.bdle.variants}
  There are yet other incarnations of the selfsame bundle that the components of $T$ are sections of, some manifesting $g$-dependence and some not. For instance, the bundle (in the usual sense \cite[Definition 37.1]{zbMATH01077733} of Fr\'echet manifolds) $\cP\xrightarrowdbl{}\cM_{\bG}(M)\times M/\bG$ restricting to
  \begin{equation*}
    \cP|_{\left\{g\right\}\times M/\bG}
    \cong
    \End\left(\ad(M)\right)_{sym}\otimes T^*\left(M/\bG\right)
  \end{equation*}
  with the ``sym'' subscript denoting symmetric matrices with respect to the transposition induced by $g$ is naturally isomorphic to the pullback
  \begin{equation*}
    \pi_2^* T\left(M/\bG\right)\otimes S^2 T^*\left(\ad(M)\right)
    ,\quad
    \cM_{\bG}(M)\times M/\bG
    \xrightarrowdbl[\text{second projection}]{\quad\pi_2\quad}
    M/\bG.
  \end{equation*}
  Such natural identifications permit unambiguous comparisons and notions of neighborhood and closeness in and between otherwise ostensibly $g$-dependent constructs such as $\End\left(\ad(M)\right)_{sym}$ and its fibers. 
\end{remark}

\begin{theorem}\label{th:gen.t.tens}
  Let $M\circlearrowleft \bG$ be a free $C^{\infty}$ action by a compact Lie group on a $C^{\infty}$ compact connected manifold and $p\in \bM$.

  For any $g\in \cM_{\bG}(M)$ and any $C^{\infty}$ neighborhood $U\ni g$, each of the two components
  \begin{equation*}
    \overrightarrow{T^{g'}_{p\bG}}
    \quad\text{and}\quad
    T_{p\bG}^{g'\uparrow}
    \quad
    \in
    \quad
    \begin{gathered}
      \End_{p\bG}(\ad(M))_{sym}\otimes T^*_{p\bG}(M/\bG)
      \cong\\
      \cong
      T_{p\bG}(M/\bG)\otimes S^2 T_{p\bG}^{*}\left(\ad(M)\right)
    \end{gathered}    
    ,\quad
    g'\in U
  \end{equation*}
  can be chosen arbitrarily within a sufficiently small neighborhood of the respective $g$ analogue being deformed. 
\end{theorem}
\begin{proof}
  The two components are interchanged by the canonical $g$-parametrized $M/\bG$-based bundle isomorphisms
  \begin{equation*}
    T\left(M/\bG\right)\otimes T^{*\otimes 2}\ad(M)
    \cong
    T^{*}\ad(M)\otimes T^*\left(M/\bG\right)\otimes T\ad(M)
  \end{equation*}
  alluded to preceding the proof by the selfsame combination \cite[(9.18c)+(9.18d)]{bss_einst_2008} cited previously. It thus suffices to handle the one (horizontal-valued) component
  \begin{equation*}
    \overrightarrow{T^{g'}}
    \in
    \Gamma
    \left(
      T(M/\bG)\otimes S^2 T^{*}\left(\ad(M)\right)
    \right).
  \end{equation*}
  The proof is very much modeled on that of \Cref{th:gen.onto.curv}, hinging on the onto-differential characterization for Fr\'echet submersions. Specifically, consider the smooth map
  \begin{equation*}
    \cM_{\bG}(M)
    \supseteq U\ni
    g'
    \xmapsto{\quad\text{Levi-Civita}\quad}
    \overrightarrow{\nabla^{g'}}|_{S^2 T\left(\ad(M)\right)}
    \in
    \Hom
    \left(
      S^2 T_{p\bG}\left(\ad(M)\right)
      ,\
      T_{p\bG}(M/\bG)
    \right)
  \end{equation*}
  Should its differential be onto at the original $g\in U$ subject to the deformation, we could once again conclude \cite[\S 1.56]{Schm23} that the map is a submersion locally around $g$ and its image (thus also completing the present proof). 
  
  The Levi-Civita first-variation formula \cite[Theorem 1.174(a)]{bss_einst_2008} on perturbing a metric gives
  \begin{equation}\label{eq:lc.var}
    \begin{aligned}
      g\left(D_{g}\overrightarrow{\nabla^{\bullet}}(\varepsilon)(X,Y),Z\right)
      &=
        \frac 12\left(\nabla^g_X\varepsilon(Y,Z)+\nabla^g_Y\varepsilon(X,Z)-\nabla^g_Z\varepsilon(X,Y)\right)\\
      \varepsilon
      &\in \Gamma\left(S^2 T^*M\right)^{\bG}
        ,\quad
        X,Y\in \Gamma\left(T^{\uparrow}M\right)
        ,\quad
        Z\in \Gamma\left(\overrightarrow{T}M\right).
    \end{aligned}   
  \end{equation}
  The requisite arbitrariness at a specific $p\in M$ is again immediate by the selfsame Whitney extension \cite[\S I.4, Theorem 4.1]{malgr_ideals} invoked once before:
  \begin{itemize}[wide]
  \item On the one hand, \Cref{eq:lc.var}'s third term is specifiable arbitrarily at $p\in M$, extensible therefrom locally around $p$ on a horizontal slice for the $\bG$-action by Whitney, and then extended $\bG$-equivariantly in a tube around $p\bG$ by $\bG$-translation. Local extensibility suffices, for the extension in question can always be chosen vanishing outside a sufficiently small $p\bG$-surrounding tube.

  \item As to the other two terms in \Cref{eq:lc.var}, their vanishing can always be ensured. 
  \end{itemize}
\end{proof}

It will also be useful to note that the two arbitrary-differential arguments employed in proving \Cref{th:gen.onto.curv,th:gen.t.tens} can be so employed simultaneously, providing the following enhancement.  

\begin{theorem}\label{th:r.and.t}
  Let $M\circlearrowleft \bG$ be a free $C^{\infty}$ action by a compact Lie group on a $C^{\infty}$ compact connected manifold and $p\in \bM$. 

  For any $g\in \cM_{\bG}(M)$ and any $C^{\infty}$ neighborhood $U\ni g$
  \begin{equation*}
    \overrightarrow{T^{g'}_{p\bG}}
    \oplus
    R^{\cH_{g'}}_{p\bG}
    \quad
    \in
    \quad
    \begin{gathered}
      T_{p\bG}(M/\bG)\otimes S^2 T_{p\bG}^{*}\left(\ad(M)\right)
      \\
      \oplus
      \\
      T_{p\bG}\left(\ad(M)\right)\otimes \Lambda^2 T^*_{p\bG}(M/\bG) 
    \end{gathered}
    ,\quad
    g'\in U
  \end{equation*}
  can be chosen arbitrarily within a sufficiently small neighborhood of the respective $g$ analogue.  \qedhere
\end{theorem}



\addcontentsline{toc}{section}{References}

\def\polhk#1{\setbox0=\hbox{#1}{\ooalign{\hidewidth
  \lower1.5ex\hbox{`}\hidewidth\crcr\unhbox0}}}


\Addresses

\end{document}